\documentclass{amsart}
\usepackage{amssymb,amstext,amsmath,amscd,amsthm,amsfonts,enumerate,graphicx,latexsym,stmaryrd,multicol,bm,xcolor}
\usepackage{hyperref}
\usepackage{cleveref}
\hypersetup{colorlinks=true}
\usepackage[all]{xy}
\usepackage{mathdots}
\usepackage{tikz}
\newtheorem{thm}{Theorem}[section]
\newtheorem{lem}[thm]{Lemma}
\newtheorem{prop}[thm]{Proposition}
\newtheorem{cor}[thm]{Corollary}
\theoremstyle{definition}
\newtheorem{dfn}[thm]{Definition}
\newtheorem{ques}[thm]{Question}
\newtheorem{rem}[thm]{Remark}
\newtheorem{eg}[thm]{Example}
\newtheorem{nota}[thm]{Notation}
\theoremstyle{remark}
\newtheorem*{ac}{Acknowledgments}
\newtheorem*{conv}{Convention}

\numberwithin{equation}{thm}
\def\C{\mathbb{C}}
\def\cok{\operatorname{Cok}}
\def\curv{\operatorname{curv}}
\def\cx{\operatorname{cx}}
\def\depth{\operatorname{depth}}
\def\edim{\operatorname{edim}}
\def\Ext{\operatorname{Ext}}
\def\Hom{\operatorname{Hom}}
\def\image{\operatorname{Im}}
\def\ker{\operatorname{Ker}}
\def\m{\mathfrak{m}}
\def\n{\mathfrak{n}}
\def\N{\mathbb{N}}
\def\ord{\operatorname{ord}}
\def\p{\mathfrak{p}}
\def\q{\mathfrak{q}}
\def\rank{\operatorname{rank}}
\def\RHom{\mathbf{R}\operatorname{Hom}}
\def\Soc{\operatorname{Soc}}
\def\Spec{\operatorname{Spec}}
\def\Tor{\operatorname{Tor}}
\def\Z{\mathbb{Z}}

\def\bsa{\boldsymbol{a}}
\def\bsb{\boldsymbol{b}}
\def\bsc{\boldsymbol{c}}
\def\bst{\boldsymbol{t}}
\def\bsx{\boldsymbol{x}}
\def\bsy{\boldsymbol{y}}

\begin{document}
\allowdisplaybreaks

\title{On some problems concerning the growth of Betti numbers}
\author{Kaito Kimura}
\address{Department of Mathematics, Purdue University, 150 N. University Street, West Lafayette, IN 47907, USA}
\email{m21018b@gmail.com}

\begin{abstract} 
We investigate several problems concerning the growth of Betti numbers over commutative noetherian rings. 
First, we establish uniform positive lower and finite upper bounds for the ratios of consecutive Betti numbers of finitely generated modules of infinite projective dimension over an arbitrary local ring, giving an affirmative answer to a question of Avramov. 
We also prove the upper semicontinuity of the curvature of residue fields of localizations over excellent rings.
Next, using zeros of Hilbert series, we obtain criteria for exponential growth of Betti and Bass numbers of graded modules. 
As applications, we establish exponential growth of the Betti numbers of canonical modules for non-Gorenstein positively graded Cohen–Macaulay domains. 
Furthermore, we show that the first Betti number of the canonical module of a non-Gorenstein generically Gorenstein Cohen–Macaulay local ring is strictly larger than the zeroth. 
In particular, both questions of Jorgensen and Leuschke have affirmative answers for Cohen–Macaulay quotients of Gorenstein local rings by monomial ideals in regular sequences.
\end{abstract}

\maketitle

\section{Introduction}

Throughout the present paper, we assume that all rings are commutative and noetherian. 
In this paper, we investigate some problems concerning the growth of Betti numbers.

We begin with a problem concerning arbitrary finitely generated modules. Avramov \cite{Av96} asked whether, for every finitely generated module $M$ of infinite projective dimension over a local ring $R$, one has 
$$
\limsup_{n\to\infty}\frac{\beta_{n+1}^R(M)}{\beta_n^R(M)}<\infty;
$$
see \cite[Problem 4.3.1]{Av96}. This problem has been extensively studied but appears to remain open; see \cite{Av96, AGP97, Ram82, Sun94, Sun98} for instance.
These studies have provided affirmative answers in several cases, for example, when $R$ is Cohen–Macaulay or (generalized) Golod, or when $M$ has finite complete intersection dimension.
Our first result gives an affirmative answer in general. Moreover, the bound can be chosen independently of the module.

\begin{thm}\label{mainA}
Let $R$ be a local ring. Then there exists a constant $C>0$ such that 
$$\limsup\limits_{n\to\infty}\frac{\beta_{n+1}^R(M)}{\beta_n^R(M)}\le C$$ for every  finitely generated $R$-module $M$ with infinite projective dimension.
\end{thm}

\noindent In the proof, we first show the existence of a uniform positive lower bound and then combine this bound with Koszul homology to obtain a uniform upper bound. 
Next, we study the curvature of the residue field of the localization at each prime ideal. 
Under mild assumptions, we show that the function on the spectrum assigning to each prime ideal the curvature of the corresponding residue field is upper semicontinuous; see Proposition \ref{curvlocus}.
This result can be viewed as a natural generalization of the openness of the regular and complete intersection loci.

Another important direction in the study of the growth of Betti numbers concerns the Betti numbers of canonical modules. Jorgensen and Leuschke \cite{JL07} posed the following two questions:
\begin{ques}\label{intro growth ques canonical}
Let $R$ be a Cohen--Macaulay local ring with canonical module $\omega_R$.
\begin{enumerate}
\item If $R$ is not Gorenstein, must the Betti numbers of the canonical module grow exponentially?
\item Does $\beta_1(\omega_R)\le\beta_0(\omega_R)$ imply that $R$ is Gorenstein?
\end{enumerate}
\end{ques}
\noindent Here, a sequence $\{\beta_i\}$ is said to have exponential growth if there exist $1<a<b$ such that $a^i<\beta_i<b^i$ for all $i\gg 0$.
According to \cite{JL07}, Jorgensen and Leuschke originally learned of a version of the first question from C. Huneke. 
These questions have affirmative answers for various classes of rings, but remain open in general; see \cite{As26, Av12, BSWY11, CCLTY26, CSV10, DGS26, He26, JL07, LM20, NJPF25} for instance.
In this paper, we approach these two questions using different methods for each.

Focusing on the roots of Hilbert series, we establish a sufficient condition for the exponential growth of Betti numbers of finitely generated graded modules over positively graded rings. In particular, by applying this result to canonical modules, we obtain affirmative answers to Question \ref{intro growth ques canonical}(1).
\begin{thm}\label{mainB}
Question \ref{intro growth ques canonical}(1) has an affirmative answer for every positively graded domain over a field.
\end{thm}
\noindent The key idea of the proof is that the failure of symmetry of the Hilbert series yields a nonzero complex root inside the unit disk at which the Hilbert series of the ring has a higher order of vanishing than that of its canonical module.
Analyzing the minimal free resolution of the canonical module, we then show that its Betti numbers grow exponentially, with a lower bound determined by the reciprocal of the absolute value of this root. Our main result concerning the latter question is as follows.
\begin{thm}\label{mainC}
Question \ref{intro growth ques canonical}(2) has an affirmative answer for every generically Gorenstein ring.
\end{thm}
\noindent In the proof, we reduce to the one-dimensional case and then study a two-periodic complex constructed from a presentation matrix of the canonical module and its adjugate matrix.

Our results provide affirmative answers to the questions for a general class of rings from a perspective different from those of previous studies. 
Furthermore, applying these ideas, we obtain affirmative answers for level rings and licci quotients; see Corollaries \ref{level} and \ref{licci0<1}.
In particular, both questions are answered affirmatively for quotients by monomial ideals.
These results refine several previous results, including \cite[Corollary 1.6 and Proposition 2.7]{JL07} and \cite[Corollary 6.7]{LM20}.
\begin{thm}\label{mainD}
Let $Q$ be a Gorenstein local ring and let $I$ be an ideal of $Q$ generated by monomials in a $Q$-regular sequence.
Suppose that $R=Q/I$ is Cohen–Macaulay. 
Then Questions \ref{intro growth ques canonical}(1) and (2) have affirmative answers for $R$.
\end{thm}
\noindent 

The paper is organized as follows. 
Section 2 contains the proof of Theorem \ref{mainA}, and Section 3 deals with the openness of curvature loci.
In Section 4, we establish criteria for exponential growth of Betti and Bass numbers in terms of Hilbert series.
Section 5 reviews known results on the growth of Betti numbers of canonical modules.
Theorem \ref{mainB}, together with applications to level rings and monomial quotients, is proved in Section 6.
Finally, Section 7 establishes Theorem \ref{mainC} and derives further consequences for linkage.
Combining the results of the last two sections yields Theorem \ref{mainD}.
We conclude the introduction by collecting the conventions used throughout the paper.

\begin{conv}
Let $(R,\m,k)$ be a local ring and let $M$ be a finitely generated $R$-module.
For each integer $n\ge 0$, we denote the $n$-th Betti number and the $n$-th Bass number of $M$ by $\beta_n^R(M)$ and $\mu_R^n(M)$, respectively.
We denote by $\Omega_R^n M$ the $n$-th syzygy of $M$ in a minimal free resolution, and by $\ell_R(M)$ the length of $M$.
The complexity and curvature of $M$ are defined by
\begin{align*}
&\cx_R(M)=\inf\{ d\ge 0 \mid \text{there exists } C>0 \text{ such that } \beta_n^R(M)\le C n^{d-1} \text{ for all } n\gg 0\}, \\
&\curv_R(M)=\limsup_{n\to\infty}\sqrt[n]{\beta_n^R(M)},
\end{align*}
respectively, where $\inf\emptyset=\infty$. When there is no risk of confusion, we omit the superscripts and subscripts.
\end{conv}

\section{The ratio of consecutive Betti numbers}

In this section, we investigate the asymptotic behavior of the ratios of consecutive Betti numbers of modules of infinite projective dimension. More precisely, we consider the following problem posed in \cite{Av96}. 

\begin{ques}\label{growth ques limsup}\cite[Problem 4.3.1]{Av96}
Let $R$ be a local ring, and let $M$ be a finitely generated $R$-module with infinite projective dimension.
Is $\limsup\limits_{n\to\infty}\frac{\beta_{n+1}^R(M)}{\beta_n^R(M)}$ always finite?
\end{ques}

\noindent This problem appears to remain open. It is also stated in \cite[Definition 2.9]{LM20} that this problem remains unsolved.
This question was known to have an affirmative answer in several cases: when $R$ is Cohen--Macaulay \cite[Corollary 4]{Ram82}, when $R$ is Golod \cite[Theorem 5.3.3]{Av96}, when $M$ has finite complete intersection dimension \cite[Theorem 8.1]{AGP97}, when $R$ is generalized Golod and $\curv_R(M)>1$ \cite[Corollary, p.89]{Sun98}, and under certain restrictions on the embedding codepth or linkage number of $R$; see \cite[Theorem 8.4]{AGP97} and \cite[Theorems 1.1 and 1.2]{Sun94}. The main result of this section provides an affirmative answer to this problem. More strongly, we show that the ratios of consecutive Betti numbers are bounded by a constant independent of the module.

\begin{thm}\label{solved}
Let $(R,\m,k)$ be a local ring. Then there exist $0<a<b$ satisfying the following condition: for every finitely generated $R$-module $M$ with infinite projective dimension, one has
$$
a\beta_n^R(M)\le \beta_{n+1}^R(M) \quad \text{and} \quad \beta_{m+1}^R(M)\le b\beta_m^R(M)
$$
for all $n\ge\depth R+1$ and $m\ge \depth R+\edim R+1$. In particular, one has
$$0<a\le \liminf\limits_{n\to\infty}\frac{\beta_{n+1}^R(M)}{\beta_n^R(M)}\le \limsup\limits_{n\to\infty}\frac{\beta_{n+1}^R(M)}{\beta_n^R(M)}\le b<\infty.
$$
\end{thm}

\begin{proof}
We consider the case $\depth R=0$.
Then $\{\Soc R \cap \m^i\}_{i\ge 0}$ forms a descending chain of $k$-vector spaces. 
By Krull's intersection theorem, $\bigcap_{i\ge 0}(\Soc R \cap \m^i)\subseteq \bigcap_{i\ge 0} \m^i=0$ holds.
So there exists an integer $s$ such that $\Soc R \cap \m^s=0$.
For $n\ge 1$ and a part of a minimal free resolution 
$$
R^{\oplus\beta_{n+1}^R(M)} \xrightarrow{d_{n+1}} R^{\oplus\beta_n^R(M)} \xrightarrow{d_n} R^{\oplus\beta_{n-1}^R(M)}
$$
of $M$, one has $\Soc (R^{\oplus\beta_n^R(M)})\subseteq \ker d_n=\Omega^{n+1}M$ since all entries of $d_n$ belong to $\m$. In particular, we have $\Soc (R^{\oplus\beta_n^R(M)})\subseteq \Soc (\Omega^{n+1}M)$, and thus 
\begin{align}\label{lseq1}
\beta_n^R(M)\ell(\Soc R)\le \ell(\Soc (\Omega^{n+1}M)).
\end{align}
Since 
$$
\Soc (\Omega^{n+1}M) \cap \m^s \Omega^{n+1}M\subseteq \Soc (R^{\oplus\beta_n^R(M)}) \cap \m^s R^{\oplus\beta_n^R(M)}=0,
$$
the composition $\Soc (\Omega^{n+1}M)\to \Omega^{n+1}M \to \Omega^{n+1}M \otimes_R R/\m^s$ is injective, which implies
\begin{align}\label{lseq2}
\ell(\Soc (\Omega^{n+1}M))\le \ell(\Omega^{n+1}M \otimes_R R/\m^s)\le \ell(R^{\oplus\beta_{n+1}^R(M)} \otimes_R R/\m^s)=\beta_{n+1}^R(M)\ell(R/\m^s).
\end{align}
Setting $a=\ell(\Soc R)/\ell(R/\m^s)$, it follows from \eqref{lseq1} and \eqref{lseq2} that $a\beta_n^R(M)\le \beta_{n+1}^R(M)$ for all $n\ge 1$.

To prove the other inequality, let $K=(0\xrightarrow{f_{e+1}} K_e \xrightarrow{f_e}  \cdots \xrightarrow{f_2}  K_1 \xrightarrow{f_1} K_0 \xrightarrow{f_0} 0)$ be the Koszul complex with respect to a system of generators of $\m$, where $e=\edim R$. Consider the exact sequences
$$
0\to \ker f_i \to K_i \to \image f_i \to 0 \quad \text{and} \quad 0\to \image f_{i+1} \to \ker f_i \to H_i(K) \to 0
$$
for $0\le i\le e$. 
Since $K_i$ is free, the former yields $\Tor_j^R(\ker f_i,M)\cong \Tor_{j+1}^R(\image f_i,M)$ for all $j>0$. On the other hand, the latter induces an exact sequence
$$
\Tor_j^R(\image f_{i+1},M) \to \Tor_j^R(\ker f_i,M) \to \Tor_j^R(H_i(K),M).
$$
Since $H_i(K)$ is a $k$-vector space, and $\ker f_i$ and $\image f_{i+1}$ are locally free on the punctured spectrum, these Tor modules have finite length for $j>0$. Hence
\begin{align*}
\ell(\Tor_{j+1}^R(\image f_i,M))
=\ell(\Tor_j^R(\ker f_i,M))
&\le \ell(\Tor_j^R(\image f_{i+1},M))
+\ell(\Tor_j^R(H_i(K),M))\\
&=\ell(\Tor_j^R(\image f_{i+1},M))
+\beta_j^R(M)(\dim_k H_i(K)).
\end{align*}
Iterating this inequality, we obtain
\begin{align*}
\ell(\Tor_{j+e}^R(\image f_1,M))
&\le \ell(\Tor_{j+e-1}^R(\image f_2,M))
+\beta_{j+e-1}^R(M)(\dim_k H_1(K)) \\
&\le \ell(\Tor_{j+e-2}^R(\image f_3,M))
+\beta_{j+e-2}^R(M)(\dim_k H_2(K))
+\beta_{j+e-1}^R(M)(\dim_k H_1(K))\\
\cdots \ \ &\le \ell(\Tor_j^R(\image f_{e+1},M))
+\sum_{h=0}^{e-1} \beta_{j+h}^R(M)(\dim_k H_{e-h}(K))
\end{align*}
for $j>0$.
Since $\image f_1=\m$ and $\image f_{e+1}=0$, it follows that
\begin{equation}\label{ineqtorbetti}
\beta_{j+e+1}^R(M)=\ell(\Tor_{j+e}^R(\image f_1,M))\le\sum_{h=0}^{e-1} \beta_{j+h}^R(M)(\dim_k H_{e-h}(K))
\end{equation}
for $j>0$. For $0\le h\le e$ and the constant $a>0$ chosen above, we have
$$
\beta_{j+h}^R(M)\le a^{-1}\beta_{j+h+1}^R(M)\le \cdots\le a^{-e+h}\beta_{j+e}^R(M).
$$
Therefore, by \eqref{ineqtorbetti}, for every $j>0$,
$$
\beta_{j+e+1}^R(M)\le\sum_{h=0}^{e-1} \beta_{j+h}^R(M)(\dim_k H_{e-h}(K))\le\Bigl(\sum_{h=0}^{e-1}a^{-e+h}(\dim_k H_{e-h}(K))\Bigr)
\beta_{j+e}^R(M).
$$
Setting $b=\sum_{h=0}^{e-1}a^{-e+h}(\dim_k H_{e-h}(K))$, we obtain the desired assertion.

Now suppose that $t=\depth R>0$. Let $\bsx=x_1,\ldots,x_t$ be an $R$-regular sequence. For every $n\ge 0$ and finitely generated $R$-module $M$ of infinite projective dimension,
$$
\beta_{n+t}^R(M)=\beta_n^R(\Omega^t M)=\beta_n^{R/\bsx R}(\Omega^t M/\bsx\Omega^t M)
$$
Consequently, the constants $a$ and $b$ obtained for $R/\bsx R$ also satisfy the desired assertion for $R$.
\end{proof}

The perspective presented in the following remark was suggested by Linquan Ma.

\begin{rem}\label{spectral sequence}
The inequality in \eqref{ineqtorbetti} can be understood more easily from the perspective of a spectral sequence.
Using the same notation as in the proof of Theorem \ref{solved}, consider 
$$
E^2_{p,q}=\Tor_p^R(M, H_q(K)) \Longrightarrow H_{p+q}(M\otimes_R^{\mathbf{L}} K)
$$
and the $E^2$-page in the $(p,q)$-plane
$$
\xymatrix@C=5mm@R=5mm{
\cdots & E^2_{j+e-3,3} & E^2_{j+e-2,3} & E^2_{j+e-1,3} & E^2_{j+e,3} & E^2_{j+e+1,3} & \cdots\\
\cdots & E^2_{j+e-3,2} & E^2_{j+e-2,2} & E^2_{j+e-1,2} & E^2_{j+e,2} & E^2_{j+e+1,2} & \cdots\\
\cdots & E^2_{j+e-3,1} & E^2_{j+e-2,1} & E^2_{j+e-1,1} & E^2_{j+e,1} & E^2_{j+e+1,1} & \cdots\\
\cdots & E^2_{j+e-3,0} & E^2_{j+e-2,0} & E^2_{j+e-1,0} & E^2_{j+e,0} & E^2_{j+e+1,0} \ar[ull] \ar[uulll] \ar[uuullll] & \cdots.
}
$$
For $j\ge 0$, one has $H_{j+e+1}(M\otimes_R^{\mathbf{L}} K)=H_{j+e+1}(M\otimes_R K)=0$.
Therefore, we must have $E^{\infty}_{j+e+1,0}=0$, and every nonzero class in $E^2_{j+e+1,0}$ must eventually be mapped nontrivially by one of the differentials extending to the upper left.
This means that $\ell(E^2_{j+e+1,0})\le \sum_{h=1}^e \ell(E^2_{j+e-h,h})$.
This is precisely the inequality \eqref{ineqtorbetti}.
\end{rem}

\section{Openness of curvature loci}

In this section, we study the upper semicontinuity of curvature. More precisely, we prove, under mild assumptions, that the locus where the residue fields of localizations have curvature less than a fixed value is open. We first establish a lemma needed for the proof.

\begin{lem}\label{opennbh}
Let $R$ be a ring, $\p$ a prime ideal of $R$, $M$ a finitely generated $R$-module, and $n\ge 0$.
Then there are $f\in R\setminus\p$ and an exact sequence $F_n\to F_{n-1} \to \cdots \to F_0 \to M_f\to 0$ where $F_0,\ldots,F_n$ are free $R_f$-modules such that its localization at $\p$ is part of a minimal free resolution of the $R_\p$-module $M_\p$.
\end{lem}

\begin{proof}
Let $R_\p^{\oplus \beta_n}\to R_\p^{\oplus \beta_{n-1}}\to \cdots \to R_\p^{\oplus \beta_0}\to M_\p \to 0$ be a minimal free resolution of $M_\p$.
Then there exists a sequence
$$
R^{\oplus \beta_n}\xrightarrow{d_n} R^{\oplus \beta_{n-1}}\xrightarrow{d_{n-1}} \cdots \xrightarrow{d_1} R^{\oplus \beta_0}\xrightarrow{d_0} M \to 0
$$
of $R$-homomorphisms (not necessarily forming a complex) whose localization at $\p$ is isomorphic to the minimal free resolution as an $R_\p$-complex.
Thus, $\image(d_{i-1}d_i)_\p=0$ for every $1\le i\le n$.
There exists $g\in R\setminus\p$ that annihilates $\image(d_{i-1}d_i)$ for every $1\le i\le n$.
Then the localization 
$$
R_g^{\oplus \beta_n}\xrightarrow{d_n} R_g^{\oplus \beta_{n-1}}\xrightarrow{d_{n-1}} \cdots \xrightarrow{d_1} R_g^{\oplus \beta_0}\xrightarrow{d_0} M_g \to 0
$$
at $g$ is a complex. Since its localization at $\p$ is exact, similarly, there exists an element $h\in R\setminus\p$ that annihilates all homologies. Setting $f=gh$, we obtain the desired exact sequence.
\end{proof}

This lemma shows that, for any fixed $n$, the first $n$ steps of a minimal free resolution at a prime $\p$ can be extended to a neighborhood of $\p$ in such a way that they remain minimal at every prime containing $\p$ in that neighborhood.
However, the neighborhood may depend on $n$, and it is much less clear whether one can find a single neighborhood that works for the entire minimal free resolution.
Nevertheless, the main result of this section, Proposition \ref{curvlocus}, shows that such a neighborhood can be found for curvature.
Before stating the proposition, we introduce some notation.
The \textit{regular locus} (resp.\! \textit{complete intersection locus}) of $R$ is the set of prime ideals $\p$ of $R$ such that $R_\p$ is regular (resp.\! a complete intersection). For $I\subseteq R$, we denote by $V(I)$ (resp.\! $D(I)$) the set of prime ideals of $R$ containing (resp.\! not containing) $I$.

\begin{prop}\label{curvlocus}
Let $R$ be a ring. Suppose that for every prime ideal $\p$ of $R$, the regular locus of $R/\p$ contains a nonempty open subset (e.g.\! $R$ is excellent).
Then 
$$
\{\p\in\Spec R \mid \curv_{R_\p} (\kappa(\p))<r \ \text{where $\kappa(\p)$ is the residue field of $R_\p$}\}
$$
is an open subset of $\Spec R$ for every $r\ge 0$.
\end{prop}

\begin{proof}
We set $U=\{\p\in\Spec R \mid \curv_{R_\p} (\kappa(\p))<r \}$.
To show that $U$ is open, we verify the two conditions in \cite[Theorem 24.2]{Mat89}. One of these conditions is immediate. Indeed, for any prime ideals $\p\subseteq\q$, we have
$$
\curv_{R_\p} (\kappa(\p))\le \curv_{R_\q} ((R/\p)_\q)\le \curv_{R_\q} (\kappa(\q))
$$
by \cite[Propositions 4.2.4(1) and 4.2.5(1)]{Av96}.
Thus, it suffices to verify the other condition: for every $\p\in U$, $U$ contains a nonempty open subset of $V(\p)$.

Let $\p\in U$. 
Let $\bsx=x_1,\ldots,x_e$ be a system of generators of $\p$.
The short exact sequence
$$
0\to \Omega_{R_\p}^{n+1} \kappa(\p)\to R_\p^{\oplus \beta_n^{R_\p}(\kappa(\p))}\to \Omega_{R_\p}^n \kappa(\p)\to 0
$$
obtained by truncating a minimal free resolution, induces a long exact sequence 
$$
\cdots \to H_i(\bsx; R_\p^{\oplus \beta_n^{R_\p}(\kappa(\p))}) \to H_i(\bsx; \Omega_{R_\p}^n \kappa(\p)) \to H_{i-1}(\bsx; \Omega_{R_\p}^{n+1} \kappa(\p)) \to \cdots
$$
in Koszul homology. For each $n\ge 0$ and $0\le i\le e$, we set
$$
c_i=\dim_{\kappa(\p)} H_i(K(\bsx;R_\p)), \quad
h_i(n)=\dim_{\kappa(\p)}H_i(K(\bsx; \Omega_{R_\p}^n\kappa(\p))).
$$
For all $1\le i\le e$, the long exact sequence yields the inequality
$$
h_i(n)\le c_i \beta_n^{R_\p}(\kappa(\p))+h_{i-1}(n+1).
$$
Since $h_0(m)=\beta_m^{R_\p}(\kappa(\p))$ and $c_0=1$, iterating this inequality over $i$ yields
$$
h_i(n)\le \sum_{j=0}^i c_{i-j} \beta_{n+j}^{R_\p}(\kappa(\p))
$$
for all $0\le i\le e$.
Fix a number $s$ such that $\curv_{R_\p} (\kappa(\p))<s<r$. The right-hand side of 
$$
h_i(n)s^{-(n+i)}\le \Bigl(\sum_{j=0}^i c_{i-j} \beta_{n+j}^{R_\p}(\kappa(\p))\Bigr)s^{-(n+i)}=\sum_{j=0}^i c_{i-j} s^{j-i}\Bigl(\beta_{n+j}^{R_\p}(\kappa(\p))s^{-(n+j)}\Bigr)
$$
tends to $0$ as $n\to\infty$.
Hence there exists $n>0$ such that $\sum_{i=0}^eh_i(n)s^{-(n+i)}<1$.
Fix such an $n$ for the remainder of the proof.

By Lemma \ref{opennbh}, there are $f\in R\setminus\p$ and an exact sequence $F_n\to F_{n-1} \to \cdots \to F_0 \to (R/\p)_f\to 0$ where $F_0,\ldots,F_n$ are free $R_f$-modules such that its localization at $\p$ is part of a minimal free resolution of the $R_\p$-module $\kappa(\p)$.
We have $X_\p=\Omega_{R_\p}^n \kappa(\p)$, where $X=\image(F_n\to F_{n-1})$.
Since $(R/\p)_f$ is an integral domain and $H_i(\bsx; X)$ is a finitely generated $(R/\p)_f$-module of rank $h_i(n)$, we can take $g\in R\setminus\p$ such that $H_i(\bsx; X)_g\cong(R/\p)^{h_i(n)}_{fg}$ for all $0\le i\le e$ by \cite[Theorem 24.1]{Mat89}.
As the regular locus of $R/\p$ contains a nonempty open subset, there is $h\in R\setminus\p$ such that $(R/\p)_h$ is regular.
Since $V(\p)\cap D(f)\cap D(g)\cap D(h)$ contains $\p$, it is a nonempty open subset of $V(\p)$. We show that this set is contained in $U$.

Let $\q\in V(\p)\cap D(f)\cap D(g)\cap D(h)$. Consider the spectral sequence
$$
E_{pq}^2=\Tor_p^{R_\q}(H_q(\bsx; X_\q), \kappa(\q)) \Longrightarrow H_{p+q}(K(\bsx; X_\q)\otimes_{R_\q}^{\mathbf{L}} \kappa(\q))
$$
Since $\bsx$ is contained in $\q$, one has 
$$
K(\bsx; X_\q)\otimes_{R_\q}^{\mathbf{L}} \kappa(\q)\simeq X_\q \otimes_{R_\q}^{\mathbf{L}}\bigl(K(\bsx; R_\q)\otimes_{R_\q}^{\mathbf{L}} \kappa(\q)\bigr)\simeq X_\q\otimes_{R_\q}^{\mathbf{L}} \left(\bigoplus_{i=0}^e \kappa(\q)^{\oplus \binom{e}{i}}[i]\right),
$$
and thus 
$$
H_m(K(\bsx; X_\q)\otimes_{R_\q}^{\mathbf{L}} \kappa(\q))=\bigoplus_{i=0}^e\Tor_{m-i}^{R_\q}(X_\q, \kappa(\q))^{\oplus \binom{e}{i}}.
$$
The condition $\q\in D(f)$ ensures that the complex $F$ can be localized at $\q$. Since $F_\p$ is minimal, every entry of each differential of $F$ belongs to $\p R_f$. As $\q$ contains $\p$, the localization $F_\q$ is also minimal. In particular, $X_\q$ is isomorphic to $\Omega^n_{R_\q}((R/\p)_\q)$.
Also, $\q\in D(f)\cap D(g)$ shows that $H_i(\bsx; X_\q)\cong(R/\p)^{h_i(n)}_\q$ for all $0\le i\le e$.
For all $m\ge 1$, we obtain 
\begin{align*}
\beta_{m+n}^{R_\q}((R/\p)_\q)=\beta_m^{R_\q}(X_\q) &\le\dim_{\kappa(\q)}\bigl(H_m(K(\bsx; X_\q)\otimes_{R_\q}^{\mathbf{L}} \kappa(\q))\bigr) \\
&\le \sum_{p+q=m} \dim_{\kappa(\q)} E_{pq}^2 =\sum_{i=0}^e h_i(n)\beta_{m-i}^{R_\q}((R/\p)_\q).
\end{align*}
Setting $a_j=\beta_j^{R_\q}((R/\p)_\q)s^{-j}$ gives
\begin{align*}
a_{m+n}\le \Bigl(\sum_{i=0}^e h_i(n)\beta_{m-i}^{R_\q}((R/\p)_\q)\Bigr)&s^{-(m+n)}
=\sum_{i=0}^e h_i(n) s^{-(n+i)} a_{m-i}\\
&\le \Bigl(\sum_{i=0}^e h_i(n) s^{-(n+i)}\Bigr) \max_{0\le i\le e}\{a_{m-i}\}\le\max_{0\le i\le e}\{a_{m-i}\}
\end{align*}
for every $m\ge 1$ by the choice of $n$.
This means that $\{a_i\}_{i\ge 0}$ is bounded, and thus $\curv_{R_\q} ((R/\p)_\q)\le s$.
Consider the spectral sequence
$$
{}'E_{pq}^2=\Tor_p^{(R/\p)_\q}\bigl(\kappa(\q), \Tor_q^{R_\q}((R/\p)_\q, \kappa(\q))\bigr) \Longrightarrow {}'H_{p+q}=\Tor_{p+q}^{R_\q}(\kappa(\q), \kappa(\q)).
$$
As $\q$ belongs to $D(h)$, $(R/\p)_\q$ is regular and we have
$$
\beta_m^{R_\q}(\kappa(\q))=\dim_{\kappa(\q)} {}'H_m\le \sum_{p+q=m} \dim_{\kappa(\q)} {}'E_{pq}^2=\sum_{i=0}^{\dim(R/\p)_\q} \beta_i^{(R/\p)_\q}(\kappa(\q)) \beta_{m-i}^{R_\q}((R/\p)_\q).
$$
Consequently, $\curv_{R_\q} (\kappa(\q))\le\curv_{R_\q} ((R/\p)_\q)\le s<r$, and thus $\q\in U$.
\end{proof}

\begin{rem}\label{notopenle}
We remark here on the distinction between the following two subsets:
\begin{align*}
&U_{<r}=\{\p\in\Spec R \mid \curv_{R_\p} (\kappa(\p))<r \ \text{where $\kappa(\p)$ is the residue field of $R_\p$}\} \\
&U_{\le r}=\{\p\in\Spec R \mid \curv_{R_\p} (\kappa(\p))\le r \ \text{where $\kappa(\p)$ is the residue field of $R_\p$}\}.
\end{align*}
By \cite[Remark 4.2.3(2)]{Av96}, $U_{<1}$ is the regular locus of $R$. 
It is known that the regular locus of a Noetherian ring need not be open; see \cite[Example 1]{Hoc73}.
This provides evidence that the assertion of Proposition \ref{curvlocus} does not hold without any assumptions.
By \cite[Corollary 8.2.2]{Av96}, $U_{\le 1}$ is the complete intersection locus of $R$.
Since Nagata's criterion holds for the complete intersection property by \cite[Theorem 3.1]{GM78}, $U_{\le 1}$ is open under the same assumptions as in Proposition \ref{curvlocus}.
In light of this, it is natural to ask whether $U_{\le r}$ is also open. 
However, ChatGPT provided an example of an excellent ring for which $U_{\le 4}$ is not open.

Let $S=K[x_1,\ldots,x_6,y]$ be a polynomial ring over a field $K$ of characteristic zero and let $I$ be an ideal of $S$ generated by
$$
\{x_i^2\}_{1\le i\le 6}, \ \{x_ix_{i+1}\}_{1\le i\le 5}, \ x_1x_3+yx_3x_6-x_4x_6, \ x_1x_4+x_3x_6+(y-2)x_4x_6.
$$
Set $R=S/I$. The excellent ring $R$ has a unique minimal prime ideal $\p=(x_1,\ldots,x_6)R$, and $R/\p=K[y]$. Thus, for each $a\in K$, the ideal $\m_a:=\p+(y-a)R$ is a maximal ideal of $R$.
Since $R$ is finitely generated over $K[y]$, $\operatorname{Ass}_{K[y]} R$ is a finite set, which means that $y-a$ is regular on $R$ for all but finitely many $a\in K$.
For an integer $a\ge 2$, the quotient $R_{\m_a}/(y-a)R_{\m_a}$ is isomorphic to the ring ``$R$'' in \cite[Example 6.2]{ACI10}, whose Poincar\'{e} series is
$$
P_k^R(1,t)=\frac{1}{(1-2t)(1-4t)-t^{a+1}(1+t)}.
$$
Since $y-a$ is part of a minimal system of generators of $\m_aR_{\m_a}$, if $y-a$ is a nonzerodivisor, the Poincar\'{e} series of $R_{\m_a}$ is obtained by multiplying that of $R_{\m_a}/(y-a)R_{\m_a}$ by $(1+t)$. Thus, it is given by
$$
\frac{1+t}{(1-2t)(1-4t)-t^{a+1}(1+t)}
$$
As $t$ increases from $0$ to $1/4$, the function $(1-2t)(1-4t)$ decreases monotonically from $1$ to $0$, while $t^{a+1}(1+t)$ increases monotonically from $0$ to $5\cdot 4^{-a-2}$. The radius of convergence $\lambda_a$ of the Poincaré series satisfies $0<\lambda_a<1/4$, and the curvature of the residue field of $R_{\m_a}$ is $\lambda_a^{-1}>4$.
Therefore $\m_a\notin U_{\le 4}$ for infinitely many $a\ge 2$.
On the other hand, as $a\to\infty$, we have $5\cdot 4^{-a-2}\to 0$, and hence $\lambda_a\to 1/4$. Consequently, for every $\varepsilon>0$, there exists $a\ge 2$ such that $\m_a\in U_{<4+\varepsilon}$.
Since $\p\subseteq\m_a$ and $U_{<4+\varepsilon}$ is stable under generalization, $\p$ belongs to $U_{<4+\varepsilon}$ for every $\varepsilon>0$. Consequently, $\p\in U_{\le 4}$.
As $R$ is a one-dimensional ring with $\p$ as its unique minimal prime, the infinitely many maximal ideals $\m_a$ in $\Spec R\setminus U_{\le 4}$ are distinct minimal elements of this set. In particular, $\Spec R\setminus U_{\le 4}$ is not closed, and hence $U_{\le 4}$ is not open.
\end{rem}

As noted in the preceding remark, the locus considered in Proposition \ref{curvlocus} can be viewed as a natural generalization of both the regular locus and the complete intersection locus. 
In a different direction, Ikeda \cite{Ikeda26} studied a refinement of the complete intersection locus by imposing bounds on the codimension. 
In \cite[Main Theorem B]{Ikeda26}, he proved that the locus of complete intersections of codimension at most $c$ is open under exactly the same assumptions as those in Proposition \ref{curvlocus}. 
That result is the version of Proposition \ref{curvlocus} obtained by replacing curvature with complexity.
On the other hand, \cite[Main Theorem A]{Ikeda26} shows that the property of being a complete intersection of codimension at most $c$ does not satisfy Nagata's criterion for any $c\ge 1$.

\section{Exponential growth of graded modules}

In this section, we study the growth of total Betti and Bass numbers of finitely generated graded modules in terms of their Hilbert series. 
More precisely, we establish criteria for exponential growth using zeros of Hilbert series and derive a necessary condition involving cyclotomic polynomials for a module to have curvature at most one. 
We first fix some notation.

\begin{nota}\label{notation}
To avoid repeatedly introducing the same notation, we fix here some notation that will be used frequently throughout this section.
Let $R=\bigoplus_{i\ge 0}R_i$ be a positively graded ring of dimension $d$ where $R_0=k$ is a field, let $\bsx=x_1,\ldots,x_d$ be a homogeneous system of parameters of $R$ with $e_i=\deg x_i$, and let $M$ be a finitely generated graded $R$-module.
We set $\delta=\min\{j>0\mid R_j\ne0\}$.
For each integer $i\ge 0$, we define the $i$th (total) Betti number and the $i$th (total) Bass number of $M$ by
$$
\beta_i^R(M):=\dim_k\Tor_i^R(k,M)
\quad\text{and}\quad
\mu_R^i(M):=\dim_k\Ext_R^i(k,M).
$$
In this case, the Hilbert series $H_M(t)$ of $M$ can be expressed as follows:
$$
H_M(t)=\frac{h_M(t)}{\prod_{i=1}^d(1-t^{e_i})}, \ \ \text{where} \ h_M(t)\in\Z[t, t^{-1}].
$$
If $M$ is maximal Cohen--Macaulay, one has $h_M(t)=H_{M/\bsx M}(t)$, and hence all coefficients of $h_M(t)$ are nonnegative integers.
 When $R$ is Cohen--Macaulay, let $\omega_R$ be a graded canonical module of $R$.
\end{nota}

\begin{rem}\label{grloc}
Let $R$ be a positively graded ring over a field, let $\m$ denote its homogeneous maximal ideal, and let $M$ be a finitely generated graded $R$-module. 
Then the total Betti numbers of $M$ arising from a minimal graded free resolution over $R$ coincide with the Betti numbers of $M_\m$ arising from a minimal free resolution over the local ring $R_\m$.
In particular, suppose that a $k$-algebra $R$ is endowed with two different positive gradings. If the corresponding homogeneous maximal ideals coincide as maximal ideals of $R$, then the graded canonical modules associated with the two gradings have the same Betti numbers.
Henceforth, we will often tacitly identify the Betti and Bass numbers of graded modules over positively graded rings over a field with those of their localizations at the homogeneous maximal ideals.
\end{rem}

We recall a result of Stanley \cite[Theorem 4.4]{St78}, which will be used not only in the proof of the main result of this section but also in our later study of the growth of the Betti numbers of canonical modules. We state it in the following form given in \cite[Corollary 4.4.6]{BH98}, which conveniently collects the properties needed for our purposes.

\begin{lem}\label{Stanley}
{\rm (See \cite[Corollary 4.4.6]{BH98}).}
Let $R$ be Cohen--Macaulay and $\omega_R$ be as in Notation \ref{notation}.
\begin{enumerate}[\rm (1)]
\item One has $H_{\omega_R}(t)=(-1)^dH_R(t^{-1})$.
In particular, if $R$ is Gorenstein, then the equality $H_R(t)=(-1)^d t^{a(R)}H_R(t^{-1})$ holds, where $a(R)$ is the $a$-invariant of $R$.
\item If $R$ is a domain and $H_R(t)=(-1)^d t^cH_R(t^{-1})$ for some $c\in\Z$, then $R$ is Gorenstein.
\end{enumerate}
\end{lem}

We show that the roots of quotients of Hilbert series of a graded module give a criterion for exponential growth of its total Betti numbers. The same argument also yields a corresponding result for Bass numbers. Before stating the result, we introduce some notation.

\begin{dfn}
For a nonzero rational function $f(t)\in\C(t)$ and $z\in\C$, let $\ord_z(f)$ denote the order of $f$ at $z$; thus, $\ord_z(f)>0$ if and only if $z$ is a zero of $f$, while $\ord_z(f)<0$ if and only if $z$ is a pole of $f$.
In these cases, $|\ord_z(f)|$ is the multiplicity of the zero or the order of the pole, respectively.
\end{dfn}

\begin{thm}\label{root and growth}
Let $R$ and $M$ be as in Notation \ref{notation}.
Let $z\in\C$ satisfy $0<|z|<1$.
\begin{enumerate}[\rm (1)]
\item If $\ord_z(H_R(t))>\ord_z(H_M(t))$, then for every $|z|<\rho<1$, one has $\beta_i(M)\ge \rho^{-\delta i}$ for $i\gg 0$.
In particular, one has $\curv_R(M)\ge |z|^{-\delta}>1$ and thus $\cx_R(M)=\infty$.
\item  If $\ord_z(H_R(t))>\ord_z(H_M(t^{-1}))$, then for every $|z|<\rho<1$, one has $\mu_R^i(M)\ge \rho^{-\delta i}$ for $i\gg0$.
\end{enumerate}
\end{thm}

\begin{proof}
(1) A shift of $M$ changes neither its total Betti numbers nor $s=\ord_z(H_R(t)/H_M(t))$. Thus, for simplicity, we may assume that $\inf\{i\mid M_i\ne 0\}=0$.
Let $F$ be the graded minimal free resolution of $M$ and let $\Omega^iM$ be the $i$th syzygy.
Since all entries of the differentials in a minimal graded free resolution lie in the homogeneous maximal ideal, $(\Omega^iM)_j=0=(F_i)_j$ for all $j<\delta i$. 
If we write $F_i=\bigoplus_{j\ge\delta i}R(-j)^{\oplus \beta_{ij}}$, then $H_{F_i}(t)=p_i(t)H_R(t)$ holds for $p_i(t)=\sum_{j\ge \delta i} \beta_{ij}t^j\in\Z[t]$, and thus $H_{\Omega^{i+1}M}=p_iH_R-H_{\Omega^iM}$ for $i\ge 0$.
Similarly, we can write $H_{\Omega^iM}(t)=\sum_{h\ge \delta i} c_{ih}t^h$ for nonnegative integers $c_{ih}$.
As $\Omega^i M$ is the quotient of $F_i$, there is a coefficientwise inequality below between the Hilbert series of them:
\begin{equation}\label{coffeq}
0\le c_{ih}\le \sum_{j+l=h} \beta_{ij}(\dim_k R_l)\ \ \text{for all} \ \ h\ge \delta i.
\end{equation}
Put $m=\ord_z(H_M(t))$. Define an operator $\Theta$ on rational functions by
$$
\Theta=t\frac{d}{dt}.
$$
Since $z\ne0$ and $\ord_z(H_R(t))>m$, one has $(\Theta^mH_M)(z)\ne0$ and $(\Theta^m(p_iH_R))(z)=0$.
It follows from $H_{\Omega^{i+1}M}=p_iH_R-H_{\Omega^iM}$ that
$$
(\Theta^m H_{\Omega^iM})(z)=-(\Theta^m H_{\Omega^{i-1}M})(z)=\cdots=(-1)^i(\Theta^mH_M)(z)\ne 0.
$$
Hence for every $i\ge0$,
$$
|(\Theta^m H_{\Omega^iM})(z)|=|(\Theta^mH_M)(z)|=:c>0
$$
Since $\Theta^m (t^h)=h^mt^h$, we have $(\Theta^m H_{\Omega^iM})(t)=\sum_{h\ge\delta i}c_{ih}h^m t^h$.
Then \eqref{coffeq} gives
\begin{align*}
c=|(\Theta^m H_{\Omega^iM})(z)|\le\sum_{h\ge\delta i}c_{ih}h^m |z|^h&\le\sum_{h\ge\delta i} \Bigl(\sum_{j+l=h} \beta_{ij}\dim_k R_l\Bigr)(j+l)^m|z|^{j+l}\\
&\le\sum_{j\ge\delta i}\beta_{ij}\sum_{l\ge0}(\dim_k R_l)(j+l)^m|z|^{j+l}.
\end{align*}
For every $|z|<\rho<1$, we take $|z|<\sigma<\rho$. There is a constant $C>0$ such that $n^m|z|^n\le C \sigma^n$ for all $n\ge0$. We have
$$
c\le C\sum_{j\ge\delta i}\beta_{ij}\sigma^j\sum_{l\ge0}(\dim_k R_l)\sigma^l \le C\sum_{j\ge\delta i}\beta_{ij}\sigma^{\delta i}\sum_{l\ge0}(\dim_k R_l)\sigma^l \le C \beta_i(M)\sigma^{\delta i} H_R(\sigma).
$$
Hence $\beta_i(M)\ge D\sigma^{-\delta i}$ for some $D>0$ and all $i\ge0$. 
Since $\sigma<\rho$, $D\sigma^{-\delta i}>\rho^{-\delta i}$ for all $i\gg0$.

(2) The ring $R$ is a quotient of a graded polynomial ring $S=k[X_1,\ldots,X_n]$.
Let $a$ be the a-invariant of $S$. 
Then $H_S(t)=(-1)^n t^aH_S(t^{-1})$ by Lemma \ref{Stanley}(1).
It follows from \cite[Theorem 1]{ABS97} that
\begin{equation}\label{Exthilb}
\sum_{i\in\Z}(-1)^i H_{\Ext_S^i(M,S)}(t)=\frac{H_M(t^{-1})H_S(t)}{H_S(t^{-1})}=(-1)^n t^aH_M(t^{-1}).
\end{equation}

In general, for a short exact sequence $0\to A\to B\to C\to 0$ of homologically bounded complexes, the additivity of Hilbert series together with the long exact sequence in homology yields
\begin{equation}\label{sescpxhilb}
\sum_{i\in\Z}(-1)^i H_{H^i(B)}(t)=\sum_{i\in\Z}(-1)^i H_{H^i(A)}(t)+\sum_{i\in\Z}(-1)^i H_{H^i(C)}(t).
\end{equation}
Let $F=(\cdots \to F^{r-1} \to F^r \to 0)$ be the graded minimal free resolution of an $R$-complex $\RHom_S(M,S)$.
For every $0\le s\le r$, we denote by $F_{\le s}=(\cdots\to F^{s-1}\to F^s\to 0)$ the hard truncation of $F$ at degree $s$.
For $1\le s\le r$, there is the natural short exact sequence $0\to F^s[-s]\to F_{\le s}\to F_{\le s-1}\to 0$, where $F^s[-s]$ denotes the complex concentrated in degree $s$ with $F^s$ in that degree.
By \eqref{sescpxhilb}, we have
$$
\sum_{i\in\Z}(-1)^i H_{H^i(F_{\le s})}(t)=\sum_{i\in\Z}(-1)^i H_{H^i(F_{\le s-1})}(t)+ (-1)^sH_{F^s}(t).
$$
In particular, we obtain 
\begin{equation}\label{Fsum}
\sum_{i\in\Z}(-1)^i H_{H^i(F_{\le r})}(t)=\sum_{i\in\Z}(-1)^i H_{H^i(F_{\le 0})}(t)+\sum_{s=1}^r(-1)^sH_{F^s}(t).
\end{equation}
By the definition of $F$, we see that $H^i(F_{\le r})=H^i(F)=\Ext_S^i(M,S)$ for all  $i\in\Z$, and $H^j(F_{\le 0})=0$ for all $j<0$.
Set $N=H^0(F_{\le 0})$.
It follows from \eqref{Exthilb} and \eqref{Fsum} that
$$
(-1)^n t^aH_M(t^{-1})=\sum_{i\in\Z}(-1)^i H_{\Ext_S^i(M,S)}(t)=H_N(t)+\sum_{s=1}^r(-1)^sH_{F^s}(t).
$$
Since $\ord_z(H_{F^s}(t))\ge \ord_z(H_R(t))>\ord_z(H_M(t^{-1}))$ for every $1\le s\le r$, one has $\ord_z(H_N(t))=\ord_z(H_M(t^{-1}))$.
Therefore, $N$ satisfies the hypothesis of (1).

On the other hand, the complex $\RHom_S(R,S)$ is a dualizing complex for $R$.
Choose an integer $u$ such that  the shift $D:=\RHom_S(R,S)[u]$ is normalized in the sense of \cite[Definition 18.2.21]{CFH24}.
Since 
$$
\RHom_R(M,D)\simeq\RHom_S(M,S)[u]
$$
and $F$ is the graded minimal free resolution of $\RHom_S(M,S)$, 
\cite[Theorem 18.2.32(b)]{CFH24} yields
$$
\mu_R^i(M)=\rank_R F^{u-i}=\beta_{i-u}^R(N)
$$
for all $i\ge u$.
Now fix $|z|<\rho<1$ and choose $|z|<\sigma<\rho$.
By (1), we have 
$$
\mu_R^i(M)=\beta_{i-u}^R(N)\ge \sigma^{-\delta (i-u)}=(\sigma^{\delta u})\sigma^{-\delta i}
$$
for all $i\gg 0$.
Since $\sigma<\rho$, one has $(\sigma^{\delta u})\sigma^{-\delta i}\ge \rho^{-\delta i}$ for all $i\gg0$.
\end{proof}

\begin{eg}
The existence of a complex number $z$ with $0<|z|<1$ as in Theorem \ref{root and growth}(1) provides a lower bound for $\curv_R(M)$, but does not necessarily determine the value of $\curv_R(M)$ itself.
Indeed, let $R=k[x,y]/(x,y)^3$.
Then $H_R(t)=1+2t+3t^2$ with respect to the standard grading, whose two zeros $z=(-1\pm i\sqrt{2})/3$ have absolute values $1/\sqrt{3}$.
Since $H_k(t)=1$,  Theorem \ref{root and growth}(1) shows $\curv_R(k)\ge \sqrt{3}>1$.
On the other hand, $R$ is Golod and the Poincar\'{e} series of $k$ is 
$$
P_k^R(t)=\frac{1+t}{1-t-3t^2};
$$
see \cite[Proposition 5.3.4]{Av96} for instance.
The radius of convergence is $(\sqrt{13}-1)/6$, and $\curv_R(k)$ coincides with its reciprocal $(\sqrt{13}+1)/2$.
In particular, $\curv_R(k)>\sqrt{3}$.
\end{eg}

As a corollary of Theorem \ref{root and growth}, we obtain a necessary condition for a module to have curvature at most one in terms of cyclotomic polynomials.
The relationship between cyclotomic polynomials and Hilbert series has also been studied in \cite{BD21}.

\begin{cor}\label{cyclotomic defect}
Let $R$ and $\bsx$ be as in Notation \ref{notation}.
Let $M$ be a finitely generated graded $R$-module with $\curv_R(M)\le 1$.
Then, up to a unit in $\Z[t,t^{-1}]$, 
$$
\frac{h_R(t)}{\operatorname{gcd}_{\Z[t,t^{-1}]}\bigl(h_R(t),h_M(t)\bigr)}
$$
is either $1$ or a product of cyclotomic polynomials.
\end{cor}

\begin{proof}
We may assume that $p(t):=h_R(t)/\operatorname{gcd}_{\Z[t,t^{-1}]}(h_R(t),h_M(t))$ is not a unit in $\Z[t,t^{-1}]$.
Since $h_R(t)\in\Z[t]$, after multiplying by suitable powers of $t$, we may assume that both $p(t)$ and $q(t)=\operatorname{gcd}_{\Z[t,t^{-1}]}(h_R(t),h_M(t))$ belong to $\Z[t]$.
Substituting $t=0$, we obtain $1=h_R(0)=p(0)q(0)$, and thus $p(0)=\pm 1$.
Since every root $z\in\C$ of $p$ satisfies $\ord_z(h_M)<\ord_z(h_R)$, Theorem \ref{root and growth}, together with the assumption $\curv_R(M)\le1$, forces $|z|\ge 1$.
Let $z_1,\ldots,z_r\in\C$ be all the roots of $p$, counted with multiplicity.
As $p(t)\in\Z[t]$, we can write $p(t)=a\prod_{i=1}^r(t-z_i)$ for some $0\ne a\in\Z$.
We have $1=|p(0)|=|a|\prod_{i=1}^r|z_i|$.
The fact that $|a|$ and each $|z_i|$ are at least $1$ implies that they must all be equal to $1$.
Thus, one of $p$ and $-p$ is a monic polynomial with integer coefficients whose roots all lie on the unit circle. It is therefore a product of cyclotomic polynomials; see \cite[Theorem 2]{Da01} for instance.
\end{proof}

\begin{eg}\label{CI}
In the setting of Notation \ref{notation}, the Hilbert series of the residue field $k$ is $H_k(t)=1$. Hence, if $\cx_R(k)<\infty$, equivalently, if $R$ is a complete intersection, then 
$$
\frac{h_R(t)}{\operatorname{gcd}\bigl(h_R(t),\prod_{i=1}^d(1-t^{e_i})\bigr)}
$$
must, up to a unit, be either $1$ or a product of cyclotomic polynomials. 
In fact, this is true. More precisely, let $R=k[x_1,\ldots,x_n]/(f_1,\ldots,f_c)$, where $f_1,\ldots,f_c$ is a homogeneous regular sequence in $k[x_1,\ldots,x_n]$.
Setting $e_i=\deg x_i$ and $d_i=\deg f_i$, by \cite[Corollary 3.3]{St78}, we have
$$
H_R(t)=\frac{\prod_{i=1}^c(1-t^{d_i})}{\prod_{i=1}^n(1-t^{e_i})}.
$$
Hence $h_R(t)$ itself is either $1$ or a product of cyclotomic polynomials. 

The converse of Corollary \ref{cyclotomic defect} does not hold. Indeed, let $R=k[x,y]/(x^2,xy,y^3)$ be an Artinian local ring. 
Since $R$ is not a complete intersection, we have $\cx_R(k)=\infty$.
As a $k$-vector space, $R=k\oplus kx\oplus ky\oplus ky^2$. 
If we set $\deg x=\deg y=1$, then the Hilbert series is $f(t)=1+2t+t^2=(1+t)^2$, which is a product of cyclotomic polynomials.
More generally, the failure of the converse of Corollary \ref{cyclotomic defect} also follows from the following general observation.
Suppose that all the coefficients are nonnegative when we write $h_R(t)=\sum_i h_i t^i$ (for example, when $R$ is Cohen–Macaulay).
For any finitely generated graded module $M$, the module $X=\bigoplus_i M(-i)^{h_i}$ satisfies $h_X(t)=h_R(t)h_M(t)$. 
Then $X$ has the same curvature as $M$ and one has $\operatorname{gcd}_{\mathbb Z[t,t^{-1}]}(h_R(t),h_X(t))=h_R(t)$.
\end{eg}


The results obtained in this section have applications to the growth problems for Betti numbers of canonical modules studied in the subsequent sections.

\section{Growth problems for Betti numbers of the canonical module}

This section is devoted to recalling several open problems concerning the growth of the Betti numbers of canonical modules. We clarify the relationships among these problems and summarize the cases in which affirmative answers are known. These problems motivate the results presented in the subsequent sections, and reviewing them will clarify what progress those results provide toward their resolution.
We consider the following questions. 
Parts (1) and (4) of Question \ref{growth ques canonical} were posed in \cite[Introduction and Question 2.6]{JL07}, whereas part (3) was (explicitly) posed in \cite[Question 1.1]{BSWY11}.

\begin{ques}\label{growth ques canonical}
Let $R$ be a non-Gorenstein Cohen--Macaulay local ring with canonical module $\omega_R$.
\begin{enumerate}
\item Is there $1<a<b$ such that $a^i<\beta_i(\omega_R)<b^i$ for all $i\gg 0$?
\item Is $\curv_R(\omega_R)>1$?
\item Is $\cx_R(\omega_R)=\infty$?
\item Is $\beta_1(\omega_R)>\beta_0(\omega_R)$?
\end{enumerate}
\end{ques}
\noindent Note that the constant $b$ in (1) always exists.
Indeed, if $R$ is a local ring and $M$ is a finitely generated $R$-module, then there is $1<b$ such that $\beta_i^R(M)<b^i$ for all $i\gg 0$ by \cite[Corollary 4.1.5]{Av96}.
Thus, (1) is equivalent to the existence of such an $a$, or equivalently, to the condition $\liminf_{n\to\infty} \sqrt[n]{\beta_n^R(\omega_R)}>1$. 
The statements in (1), (2), and (3) become successively weaker. The converse is not known.
In fact, it is an open problem whether $\cx_R M=\infty$ implies $\curv_R M>1$ for a finitely generated $R$-module $M$; see \cite[Problem 4.3.7]{Av96}.
All of these questions remain open, and no counterexamples are known.
In contrast, Question \ref{growth ques canonical}(4) concerns the initial Betti numbers, whereas (1)--(3) concern their asymptotic growth.
It is known that $\beta_1^R(\omega_R)\ge\beta_0^R(\omega_R)$ always holds; see \cite[Lemma A.1]{CSV10}. Thus, (4) asks whether equality can occur for a non-Gorenstein Cohen--Macaulay local ring.

Several affirmative answers to Question \ref{growth ques canonical}(1) are known. 
Jorgensen and Leuschke \cite[Proposition 1.1]{JL07} proved this when $R$ is Golod, has codimension at most $3$, is one link from a complete intersection, or is of radical cube zero. 
They also obtained a general criterion for exponential growth and applied it to certain monomial quotient rings \cite[Theorem 1.5 and Corollary 1.6]{JL07}. 
In the same paper, they proved the strict inequality in (4) for Cohen--Macaulay local rings of codimension $2$ \cite[Proposition 2.7]{JL07}.
Christensen, Striuli, and Veliche \cite[Main Theorem (1.2)]{CSV10} established affirmative answers to both (1) and (4) of Question \ref{growth ques canonical} for several classes of rings, including Golod rings, nontrivial fiber products, rings with radical cube zero, Teter rings, and Artinian rings whose socles are not contained in the square of their maximal ideals. 
Avramov \cite[Theorem 4.1]{Av12} obtained the same conclusions for rings of embedding codepth at most $3$.
Borna, Sather-Wagstaff, and Yassemi \cite[Theorem 1.3]{BSWY11} showed the growth assertion in (1) for rings that are homologically of minimal multiplicity of type $(m,n,t)$ with $n>m$, and for non-Gorenstein rings that are strongly homologically of minimal multiplicity.

Another approach arises from the theory of dominant and quasi-dominant local rings developed by Takahashi \cite{T23}. 
For a non-Gorenstein Cohen--Macaulay local ring with a canonical module, quasi-dominance implies that the residue field belongs to the thick subcategory generated by the ring and its canonical module. 
Consequently, the growth condition in Question \ref{growth ques canonical}(1) holds for quasi-dominant rings.
In particular, this applies to Cohen--Macaulay Burch rings, rings with quasi-decomposable maximal ideals, excellent rings of finite Cohen--Macaulay representation type, and almost Gorenstein rings with infinite residue fields; see \cite[Propositions 5.10 and 10.15]{T23}.
More recently, Celikbas, Celikbas, Lyle, Takahashi, and Yao \cite[Corollary 4.17(ii)]{CCLTY26} proved the same growth property under the numerical condition $e(R)<2r(R)$. 
Related results were obtained by Dey, Ghosh, and Saha \cite[Corollary 5.5]{DGS26}, whose estimates for complexity and curvature yield affirmative answers to Question \ref{growth ques canonical}(2) under the weaker numerical condition $e(R)\le 2r(R)$.

Turning to Question \ref{growth ques canonical}(4), in addition to the cases discussed above, Lyle and Monta\~{n}o \cite[Corollary 6.7]{LM20} verified the desired strict inequality for generically Gorenstein rings of type $2$. 
Asgharzadeh \cite[Proposition 3.6]{As26} proved the same inequality for quotients of polynomial rings by nonprincipal Cohen--Macaulay homogeneous ideals having a $q$-linear resolution for some $q>1$. 
Nascimento, Jorge-P\'{e}rez, and Freitas \cite[Theorem 5.2]{NJPF25} settled the question for Cohen--Macaulay idealization rings, while He \cite[Corollary 3.4]{He26} established an affirmative answer for Artinian local rings satisfying $\Soc R=\m^2$ and $\edim R\ge 4$.

\begin{rem}\label{rmkques}
The conditions in Question \ref{growth ques canonical} behave well under several operations; see \cite[Remark 1.3]{JL07}.
Let $(R,\m)$ be a Cohen--Macaulay local ring with canonical module.
For a regular element $x\in\m$, $\beta_i^R(\omega_R)=\beta_i^{R/xR}(\omega_{R/xR})$ for all $i\ge 0$; hence each of (1)--(4) holds for $R$ if and only if it holds for $R/xR$.
For a flat local extension of $R\to S$ such that $S/\m S$ is Gorenstein, $\beta_i^R(\omega_R)=\beta_i^S(\omega_S)$ for all $i\ge 0$; hence each of (1)--(4) holds for $R$ if and only if it holds for $S$.
For a local ring homomorphism $R\to S$ of finite flat dimension, $\beta_i^R(\omega_R)\le \beta_i^S(\omega_S)$ for all $i\gg 0$.
Consequently, if any of (1)--(3) holds for $R$, then the same condition holds for $S$.
For example, by first passing to an infinite residue field and then taking a quotient by a reduction, the case of radical cube zero can be extended to the case of almost minimal multiplicity.
\end{rem}

By Remark \ref{grloc}, the above questions can also be posed for positively graded rings over a field. In the next section, we give an affirmative answer to Question \ref{growth ques canonical}(1), and hence to all of (1)--(3), under certain assumptions. In the following section, we investigate the inequality in (4).

\section{Exponential growth of canonical modules}

In this section, we establish criteria for exponential growth of the Betti numbers of canonical modules in terms of Hilbert series. 
We first show that the failure of symmetry of the Hilbert series also forces exponential growth of Bass numbers. 
We then apply these criteria to non-Gorenstein positively graded Cohen--Macaulay domains, level rings, and monomial quotient rings. 
The following lemma plays a central role in this section.
For $f(t)\in\Z[t]$, we define $f^\ast(t):=t^{\deg f}f(t^{-1})\in\Z[t]$.

\begin{lem}\label{small root}
Let $f(t)\in\Z[t]$ such that $f(0)=1$ and $f(1)\ne 0$.
If $f\ne f^\ast$, then there exists $z\in\C$ such that $0<|z|<1$ and $\ord_z(f)>\ord_z(f^\ast)$.
\end{lem}

\begin{proof}
Put $f(t)=\sum_{i=0}^s f_i t^i$, where $f_0,\ldots,f_s$ are integers and $s=\deg f$. Then $f^\ast=\sum_{i=0}^s f_{s-i} t^i$. 
By assumption, $f_0=f(0)=1$, and thus $f$ and $f^\ast$ are primitive. 
By Gauss's lemma, we may choose $g,p,q\in\Z[t]$ to be primitive such that $g$ is the greatest common divisor of $f$ and $f^\ast$, $f=gp$ and $f^\ast=gq$.
Since $f_0=1$, $f^\ast$ is monic. So we may assume that $g$ and $q$ are also monic, and then the leading coefficient of $p$ is $f_s$.
In what follows, we show that $p$ has a root $z\in\C$ satisfying $0<|z|<1$.

Put $a=\deg g$ and $b=s-a=\deg p=\deg q$.
Since 
$$f^\ast=t^sf(t^{-1})=(t^ag(t^{-1}))(t^bp(t^{-1}))=g^\ast p^\ast \ \text{and} \ f=(f^{\ast})^\ast=t^sf^\ast(t^{-1})=(t^ag(t^{-1}))(t^bq(t^{-1}))=g^\ast q^\ast,
$$ 
$g^\ast$ is a common divisor of $f$ and $f^\ast$.
The fact that $g$ and $g^\ast$ have the same degree and are primitive implies that $g^\ast=\pm g$.
Evaluating at $t=1$ gives the sum of the coefficients, and hence $g(1)=g^\ast(1)$. 
Since $0\ne f(1)=g(1)p(1)$, we get $g(1)\ne 0$ and $p(1)\ne 0$. Therefore, $g^\ast=g$.
One has $g(0)=1$ as $g^\ast=g$ is monic, and thus $p(0)=f(0)/g(0)=1$.
Let $z_1,\ldots,z_b\in\C$ be the roots of $p$.
The leading coefficient of $p$ is $f_s$, so $p(t)=f_s \prod_{i=1}^b(t-z_i)$ in $\C[t]$.
Evaluating at $t=0$ gives $1=|p(0)|=|f_s| \prod_{i=1}^b |z_i|$. 
Assume that $|z_i|\ge 1$ for all $1\le i\le b$.
Then $|f_s|=|z_i|=1$ for all $1\le i\le b$ as $f_s$ is an integer.
So the reciprocal of each root of $p$ is its complex conjugate, and since $p$ has real coefficients, it is also a root of $p$ with the same multiplicity.
Thus, $p$ and $p^\ast$ have the same roots, counted with multiplicity, and hence differ only by a nonzero scalar multiple.
Together with the facts that $0\ne p(1)=p^\ast(1)$, this forces $p=p^\ast$. 
But then $f=gp=g^\ast p^\ast=f^\ast$, a contradiction.
It follows that $p$ has a root $z\in\C$ with $0<|z|<1$.

Since $z$ is a root of $p$ and $g$ is the greatest common divisor of $f$ and $f^\ast$, we have 
$$
\ord_z(f)>\ord_z(f)-\ord_z(p)=\ord_z(g)=\min\{\ord_z(f),\ord_z(f^\ast)\},
$$
which means that $\ord_z(f)>\ord_z(f^\ast)$. The proof is now completed.
\end{proof}

\begin{rem}
The idea behind the proof of Lemma \ref{small root} is simple. Indeed, by factoring $f$ into linear factors and setting $t=0$, we see that the product of the absolute values of its roots is at most $1$. As we saw in the proof, the conditions $f\ne f^\ast$ and $f(1)\ne 0$ imply that $f$ must have a root of absolute value strictly less than $1$. If such a root is not a root of $f^\ast$, there is nothing more to prove. To handle the case where it is, we take the greatest common divisor of $f$ and $f^\ast$ and apply the same argument to the polynomial $p$ introduced in the proof.
\end{rem}

This lemma allows us to apply Theorem \ref{root and growth} to the problems of exponential growth of the Bass numbers of rings and the Betti numbers of canonical modules.

\begin{thm}\label{Bass growth}
Let $R$ be as in Notation \ref{notation}, and set $d=\dim R$. Suppose that $H_R(t)\ne(-1)^dt^cH_R(t^{-1})$ for every $c\in\Z$.  
Then there exists $a>1$ such that $\mu_R^i(R)>a^i$ for all $i\gg0$.
\end{thm}

\begin{proof}
With the notation of Notation \ref{notation}, 
$$
H_R(t)=\frac{h_R(t)}{\prod_{i=1}^d(1-t^{e_i})} \ \ \text{and} \ \ H_R(t^{-1})=\frac{h_R(t^{-1})}{\prod_{i=1}^d(1-t^{-e_i})}=\frac{t^{-\deg (h_R)}h_R^\ast(t)}{(-1)^dt^{-\sum_{j=1}^d e_j}\prod_{i=1}^d(1-t^{e_i})}.
$$
The assumption implies $h_R\ne h_R^\ast$. 
Moreover, since $H_R(0)=1$ and $\ord_1(H_R(t))=-d$, the polynomial $h_R$ always satisfies $h_R(0)=1$ and $h_R(1)\ne 0$, regardless of the powers of $t$ appearing in the denominator.
The assertion follows from Theorem \ref{root and growth}(2) and Lemma \ref{small root}.
\end{proof}

\begin{thm}\label{canonical growth}
Let $R$ be Cohen--Macaulay and $\omega_R$ be as in Notation \ref{notation}.
If $H_R(t)\ne (-1)^d t^{a(R)}H_R(t^{-1})$, then there exists $a>1$ such that
$\beta_i(\omega_R)>a^i$ for all $i\gg 0$. 
\end{thm}

\begin{proof}
Set $W=\omega_R(-a(R))$; that is, we shift $\omega_R$ by the $a$-invariant so that $\inf\{i\in\Z\mid W_i\ne 0\}=0$.
It follows from Lemma \ref{Stanley} that
$$
\frac{H_{W/\bsx W}(t)}{\prod_{i=1}^d(1-t^{e_i})}=H_W(t)=t^{a(R)}H_{\omega_R}(t)=t^{a(R)}(-1)^dH_R(t^{-1})=\frac{t^{e+a(R)} h_R(t^{-1})}{\prod_{i=1}^d(1-t^{e_i})},
$$
where $e=\sum_{i=1}^de_i$.
Comparing the lowest degrees of the numerator polynomials, we must have
$e+a(R)=\deg h_R(t)$, and thus $h_W(t)=h_R^\ast(t)$. The assumption $H_R(t)\ne (-1)^d t^{a(R)}H_R(t^{-1})=H_W(t)$ yields that $h_R\ne h_W=h_R^\ast$.
The assertion follows from Theorem \ref{root and growth}(1) and Lemma \ref{small root}.
\end{proof}

The following is an immediate consequence of Lemma \ref{Stanley}(2) and Theorem \ref{canonical growth}.

\begin{cor}\label{graded domain} 
The Betti numbers of the canonical module of a non-Gorenstein positively graded Cohen--Macaulay domain over a field grow exponentially.
\end{cor}

The ring $R$ is called \textit{level} if all elements in a minimal set of generators of $\omega_R$ have the same degree.

\begin{cor}\label{level}
If $R$ is Cohen--Macaulay and $h_R(t)$ in Notation \ref{notation} is not monic, then there exists $a>1$ such that
$\beta_i(\omega_R)>a^i$ for all $i\gg 0$.
In particular, the Betti numbers of the canonical module of a non-Gorenstein positively graded Cohen--Macaulay level ring over a field grow exponentially.
\end{cor}

\begin{proof}
Since $h_R(0)=1$, one has $h_R(t)\ne h_R^\ast(t)$ if $h_R(t)$ is not monic.
The first assertion follows from Theorem \ref{root and growth}(1) and Lemma \ref{small root}.
If $R$ is level, then so is $R/\bsx R$.
When we write $h_R(t)=H_{R/\bsx R}(t)=\sum_{i=0}^s h_it^i$ with $h_s\ne 0$, one has $h_s=\dim_k \Soc(R/\bsx R)>1$ since $R$ is level but not Gorenstein.
Hence $h_R(t)$ is not monic. The conclusion follows from the first assertion.
\end{proof}

Along with Corollary \ref{graded domain}, the following proposition is one of the most important applications of Theorem \ref{canonical growth}. 
The key idea of the proof is to find a grading for which the Hilbert series fails to be symmetric by assigning different weights to the variables.

\begin{prop}\label{monomial}
Let $k[x_1,\ldots,x_n]$ be a polynomial ring over $k$ and let $I$ be a monomial ideal such that $R=k[x_1,\ldots,x_n]/I$ is a non-Gorenstein Cohen--Macaulay ring. Then the Betti numbers of the canonical module $\omega_R$ of $R$ grow exponentially.
\end{prop}

\begin{proof}
Any variable that does not appear in the generators of $I$ is a regular element on $R$, and hence can be eliminated by passing to the quotient, by Remarks \ref{rmkques} and \ref{grloc}. Therefore, we may assume that every variable appears in some generator of $I$.
By \cite[Proposition 1.6.1]{HH11}, there exists a polynomial ring $T$ obtained from $k[x_1,\ldots,x_n]$ by adjoining variables and a squarefree monomial ideal $J$ of $T$, called a polarization of $I$, such that $(T/J)/\bsy(T/J)\cong R$ for some $T/J$-regular sequence $\bsy$ of linear forms.
Again, by Remarks \ref{rmkques} and \ref{grloc}, replacing $R$ by $T/J$, we may assume that $I$ is squarefree.

In this case, there exists a simplicial complex $\Delta$ such that $I=I_{\Delta}$ and $R=k[\Delta]$. 
Let $\boldsymbol{d}_1,\ldots,\boldsymbol{d}_n$ denote the standard basis of $\Z^n$, and regard $k[\Delta]$ as an $\N^n$-graded ring by setting $\deg x_i=\boldsymbol{d}_i$. 
Let $\bst=t_1,\ldots,t_n$ be indeterminates. 
The Hilbert series of $k[\Delta]$ is defined by 
$$
H_{k[\Delta]}(\bst)=H_{k[\Delta]}(t_1,\ldots,t_n)=\sum_{\bsa\in\Z^n}
\dim_k k[\Delta]_{\bsa}\,
\bst^{\bsa}=\sum_{F\in\Delta} \prod_{i\in F} \frac{t_i}{1-t_i};
$$
see the discussion following \cite[Example 5.1.6]{BH98}.
Since $R$ is a non-Gorenstein Cohen--Macaulay ring, it follows from \cite[Exercise 5.6.8]{BH98} that 
\begin{equation}\label{monom1}
H_{k[\Delta]}(\bst^{-1})=H_{k[\Delta]}(t_1^{-1},\ldots,t_n^{-1})\ne(-1)^d\bst^{\bsa}H_{k[\Delta]}(\bst)
\end{equation}
for any $\bsa=(a_1,\ldots,a_n)\in\Z^n$. 
When we set 
\begin{align}
&K(\bst):=H_{k[\Delta]}(\bst)\prod_{j=1}^n (1-t_j)=\sum_{F\in\Delta} \biggl(\prod_{i\in F} t_i\biggr)\biggl(\prod_{j\notin F} (1-t_j)\biggr) \in\Z[\bst], \label{K(t)} \\
&L(\bst):={(-1)^{d+n}}(t_1\cdots t_n)K(\bst^{-1})=(-1)^{d+n}\sum_{F\in\Delta}\biggl(\prod_{j\notin F} (t_j-1)\biggr)\in\Z[\bst]. \label{L(t)}
\end{align}
one has
\begin{equation}\label{monom2}
H_{k[\Delta]}(\bst)=\frac{K(\bst)}{\prod_{j=1}^n (1-t_j)}, \ \
H_{k[\Delta]}(\bst^{-1})=\frac{K(\bst^{-1})}{\prod_{j=1}^n (1-t_j^{-1})}=\frac{(-1)^dL(\bst)}{\prod_{j=1}^n (1-t_j)}.
\end{equation}
By \eqref{monom1}, for any $\bsa=(a_1,\ldots,a_n)\in\Z^n$, 
\begin{equation}\label{monom3}
\bst^{\bsa}K(\bst)\ne L(\bst).
\end{equation}

We consider the $\N$-graded structure on $k[\Delta]$ obtained by assigning the weight $\deg x_i=3^{i-1}$ to each variable $x_i$.
Let $\widehat{H}(t)$ denote the Hilbert series of the $\N$-graded ring $R=k[\Delta]$ with respect to this grading. 
Then, by substituting $\bst=(t^{3^0},t^{3^1},\ldots,t^{3^{n-1}})$, we obtain 
\begin{align*}
&\widehat{H}(t)=H_{k[\Delta]}(t,t^3,\ldots,t^{3^{n-1}})=\frac{K(t,t^3,\ldots,t^{3^{n-1}})}{\prod_{j=1}^n (1-t^{3^{j-1}})} \\
&\widehat{H}(t^{-1})=H_{k[\Delta]}(t^{-1},t^{-3},\ldots,t^{-3^{n-1}})=\frac{(-1)^dL(t,t^3,\ldots,t^{3^{n-1}})}{\prod_{j=1}^n (1-t^{3^{j-1}})}, 
\end{align*}
by \eqref{monom2}.
By Theorem \ref{canonical growth}, it is enough to show that $\widehat{H}(t)\ne (-1)^d t^a\widehat{H}(t^{-1})$ for all $a\in\Z$.
Suppose otherwise. 
Then there exists $a\in\Z$ such that $K(t,t^3,\ldots,t^{3^{n-1}})=t^aL(t,t^3,\ldots,t^{3^{n-1}})$.
The map
$$
\Phi: \{0,1,2\}^n\to\N: (a_1,\ldots,a_n)\mapsto\sum_{i=1}^na_i3^{i-1}
$$
is injective. 
Indeed, this is simply part of the correspondence given by ternary representation.
By \eqref{K(t)} and \eqref{L(t)}, every monomial appearing with nonzero coefficient in $K(\bst)$ or $L(\bst)$ is squarefree. 
Hence, by the injectivity of $\Phi$, the distinct monomials occurring in $K(\bst)$ and $L(\bst)$ are sent to distinct monomials.
By the equality $K(t,t^3,\ldots,t^{3^{n-1}})=t^aL(t,t^3,\ldots,t^{3^{n-1}})$, the monomials of $K(\bst)$ and $L(\bst)$ are in one-to-one correspondence.
If $\bst^{\bsb}$ and $\bst^{\bsc}$ are corresponding monomials of $K(\bst)$ and $L(\bst)$, then
$$
a=\sum_{i=1}^n(b_i-c_i)3^{i-1}, \ \text{where} \ \bsb=(b_1,\ldots,b_n), \ \bsc=(c_1,\ldots,c_n),
$$
and therefore
$$
a+\sum_{i=1}^n3^{i-1}=\sum_{i=1}^n(b_i-c_i+1)3^{i-1}.
$$
Since $b_i,c_i\in\{0,1\}$, we have $b_i-c_i+1\in\{0,1,2\}$, so the injectivity of $\Phi$ shows that $\bsb-\bsc$ is independent of the chosen pair. 
Thus there exists $\bsa\in\Z^n$ such that $\bst^{\bsa}K(\bst)=L(\bst)$, contrary to \eqref{monom3}. 
\end{proof}

We provide an example in which the symmetry of the Hilbert series changes depending on the choice of grading.
As observed in Remark \ref{grloc}, if there exists even one grading for which the Hilbert series fails to be symmetric, then the Betti numbers of the canonical module grow exponentially.

\begin{eg}
(1) The ring $R=k[x,y]/(x^2,xy,y^3)$ considered in Example \ref{CI} has Hilbert series $f(t)=1+2t+t^2$ with respect to the grading $\deg x=1=\deg y$, and this satisfies $f=f^\ast$.
On the other hand, if we set $\deg x=1$ and $\deg y=2$, then the Hilbert series is
$g(t)=1+t+t^2+t^4$, and in this case $g^\ast(t)=1+t^2+t^3+t^4\ne g(t)$.

(2) As an example showing that the domain assumption in Lemma \ref{Stanley}(2) cannot even be weakened to the assumption that the ring is reduced, \cite[Example 4.5]{St78} considers the quotient ring $R=k[x,y,z,w]/(xw,yw,xyz)$ by a squarefree monomial ideal.
This is a two-dimensional reduced Cohen--Macaulay ring, and, with the standard grading, its Hilbert series is 
$$
\frac{1+2t+t^2}{(1-t)^2},
$$
which is symmetric. 
However, if we adopt the grading $\deg x=\deg y=\deg z=1$ and $\deg w=2$, then its Hilbert series is 
$$
\frac{1+2t+3t^2+t^3}{(1-t)(1-t^2)},
$$
which is not symmetric.
Indeed, $x-z, w-x^2-y^2$ form an $R$-regular sequence of degrees $1$ and $2$, respectively, and 
$$
R/(x-z,w-x^2-y^2)\cong k[x,y]/(x^3+xy^2,y^3,x^2y).
$$
The latter ring has Hilbert series $1+2t+3t^2+t^3$, as is easily seen. 
\end{eg}

Proposition \ref{monomial} can be stated in a more general form. Indeed, for a quotient of any Gorenstein ring by a monomial ideal in a regular sequence, the Betti numbers of its canonical module depend only on the combinatorial structure of the monomial ideal.

\begin{lem}\label{samemonosameBetti}
Let $(Q,\n,k)$ be a Gorenstein local ring, let $x_1,\ldots,x_n$ be a $Q$-regular sequence, and let $I$ be an ideal of $Q$ generated by monomials in $x_1,\ldots,x_n$.
Let $k[X_1,\ldots,X_n]$ be a polynomial ring over $k$ and let $J$ be the monomial ideal corresponding to $I$.
Set $R=Q/I$ and $S=k[X_1,\ldots,X_n]/J$.
Then $R$ is Cohen--Macaulay (resp.\! Gorenstein) if and only if so is $S$, and then $\beta_i^R(\omega_R)=\beta_i^S(\omega_S)$ for all $i\ge 0$.
\end{lem}

\begin{proof}
Set $A=(Q[X_1,\ldots,X_n]/J Q[X_1,\ldots,X_n])_{(\n,X_1,\ldots,X_n)}$.
Then $X_1-x_1,\ldots,X_n-x_n$ is an $A$-regular sequence with quotient $R$; see the proof of \cite[Corollary 4.6]{GK26} for instance.
Hence $R$ is Cohen--Macaulay (resp.\! Gorenstein) if and only if so is $A$.
The latter is equivalent to $S$ being Cohen--Macaulay (resp.\! Gorenstein) since $Q\to A$ is flat with closed fiber $S_{(X_1,\ldots,X_n)}$.
In the Cohen--Macaulay case, $\beta_i^R(\omega_R)=\beta_i^A(\omega_A)$ for all $i\ge 0$ by Remark \ref{rmkques}.
As $Q$ is Gorenstein, applying \cite[Theorem 17.4.15]{CFH24} to $Q\to A$ yields the following identities for the Poincar\'{e} and Bass series:
$$
t^{\dim A}P_{\omega_A}^A(t)=I_A(t)=I_Q(t)I_S(t)=(t^{\dim Q})(t^{\dim S}P_{\omega_S}^S(t))=t^{\dim Q+\dim S}P_{\omega_S}^S(t).
$$
Since $\dim A=\dim Q+\dim S$, we have $\beta_i^A(\omega_A)=\beta_i^S(\omega_S)$ for all $i\ge 0$.
\end{proof}

As a consequence of Proposition \ref{monomial} and Lemma \ref{samemonosameBetti}, we obtain the following result.
We weaken the assumption that the base ring is regular in \cite[Corollary 1.6]{JL07} to the assumption that it is Gorenstein, and extend the result to arbitrary monomial ideals.

\begin{cor}\label{regseqmonom}
Let $(Q,\n,k)$ be a Gorenstein local ring, let $x_1,\ldots,x_n$ be a $Q$-regular sequence, and let $I$ be an ideal of $Q$ generated by monomials in $x_1,\ldots,x_n$.
If $R=Q/I$ is a non-Gorenstein Cohen--Macaulay ring, then the Betti numbers of the canonical module $\omega_R$ of $R$ grow exponentially.
\end{cor}

\section{The first two Betti numbers}

In this final section, we turn from asymptotic growth to the first two Betti numbers of canonical modules. Our main theorem rules out equality between them for non-Gorenstein generically Gorenstein Cohen--Macaulay local rings, without any further assumptions. We then obtain consequences for monomial quotient rings and ideals in the linkage class of a generically complete intersection.

\begin{thm}\label{genGor0<1}
Let $R$ be a Cohen--Macaulay local ring with canonical module $\omega_R$.
If $R$ is not Gorenstein but generically Gorenstein, then $\beta_1(\omega_R)>\beta_0(\omega_R)$.
\end{thm}

\begin{proof}
The non-Gorenstein locus of $R$ (equivalently, the non-free locus of $\omega_R$) is closed. 
Since $R$ is generically Gorenstein, every height one prime belonging to this locus is minimal in the locus. 
As the locus is closed, there are only finitely many such height one primes.
If $\dim R>1$, we can choose an $R$-regular element that is not contained in any height one prime belonging to the non-Gorenstein locus of $R$.
Then $R/xR$ is again generically Gorenstein.
Since Cohen–Macaulayness, the non-Gorenstein property, and the Betti numbers of the canonical module are preserved, we may replace $R$ by $R/xR$. 
Repeating this process, we may assume that $\dim R=1$.
Since the property that the canonical module is locally free on the punctured spectrum is preserved under completion, we may also assume that $R$ is complete.

By \cite[Proposition 2.3]{JL07} (or \cite[Lemma A.1]{CSV10}), it is known that $\beta_1(\omega_R)\ge \beta_0(\omega_R)$. 
Suppose that equality holds. We derive a contradiction.
Put $t=\beta_1(\omega_R)=\beta_0(\omega_R)$. There exists an exact sequence 
$$
R^{\oplus t}\xrightarrow{A} R^{\oplus t}\to \omega_R\to 0.
$$
Let $B$ be the adjugate of $A$. Then $AB=(\det A)I_t=BA$.
Since $(\det A)R^{\oplus t}=\image (AB)\subseteq \image A$, $\det A$ annihilates $\cok A=\omega_R$, and thus $\det A=0$.
This implies that the following is a complex:
\begin{equation}\label{ABAB}
\cdots R^{\oplus t}\xrightarrow{A} R^{\oplus t}\xrightarrow{B}R^{\oplus t}\xrightarrow{A} R^{\oplus t}\xrightarrow{B} R^{\oplus t} \cdots.
\end{equation}
For every non-maximal prime ideal $\p$ of $R$, $R_\p$ is Gorenstein as $R$ is one-dimensional generically Gorenstein Cohen--Macaulay local ring.
Since $(\omega_R)_\p$ is a free module of rank one, the matrix $A_\p$ is equivalent to the diagonal matrix $\operatorname{diag}(1,\ldots,1,0)$ of rank $t-1$.
Consequently, $B_\p$ is equivalent to the rank-one diagonal matrix $\operatorname{diag}(0,\ldots,0,1)$, and the localization of \eqref{ABAB} at $\p$ is split acyclic.

For each finitely generated $R$-module $M$, we set 
$$
X_M=H(M^{\oplus t}\xrightarrow{B_M} M^{\oplus t}\xrightarrow{A_M}M^{\oplus t})\quad \text{and} \quad Y_M=H(M^{\oplus t}\xrightarrow{A_M} M^{\oplus t}\xrightarrow{B_M}M^{\oplus t}),
$$
where $A_M$ and $B_M$ denote the endomorphisms of $M^{\oplus t}$ induced by the actions of the matrices $A$ and $B$, respectively.
Since the localization of \eqref{ABAB} at every nonmaximal prime ideal is split acyclic, both $X_M$ and $Y_M$ have finite length.
We set 
$$
\theta(M):=\ell_R(X_M)-\ell_R(Y_M).
$$
We claim that $\theta(M)=0$ for every finitely generated $R$-module $M$.
If $M$ has finite length, then 
\begin{align*}
\ell_R(X_M)=\ell_R(\ker A_M)-\ell(\image B_M)&=\bigl(\ell_R(M^{\oplus t})-\ell_R(\image A_M)\bigr)-\bigl(\ell_R(M^{\oplus t})-\ell_R(\ker B_M)\bigr)\\
&=\ell_R(\ker B_M)-\ell(\image A_M)=\ell_R(Y_M),
\end{align*}
and hence $\theta(M)=0$.
A short exact sequence $0\to L\to M\to N\to 0$ induces the short exact sequence
\[
  \xymatrix@C=20pt@R=15pt
  {
    \cdots \ar[r]^{B_L}
    & L^{\oplus t} \ar[r]^{A_L} \ar[d]
    & L^{\oplus t} \ar[r]^{B_L} \ar[d]
    & L^{\oplus t} \ar[r]^{A_L} \ar[d]
    & L^{\oplus t} \ar[r]^{B_L} \ar[d]
    & \cdots \\
    \cdots \ar[r]^{B_M}
    & M^{\oplus t} \ar[r]^{A_M} \ar[d]
    & M^{\oplus t} \ar[r]^{B_M} \ar[d]
    & M^{\oplus t} \ar[r]^{A_M} \ar[d]
    & M^{\oplus t} \ar[r]^{B_M} \ar[d]
    & \cdots \\
    \cdots \ar[r]^{B_N}
    & N^{\oplus t }\ar[r]^{A_N} 
    & N^{\oplus t} \ar[r]^{B_N}
    & N^{\oplus t} \ar[r]^{A_N}
    & N^{\oplus t} \ar[r]^{B_N}
    & \cdots \\
  }
\]
of complexes, which induces the long exact sequence 
$$
\cdots Y_N \xrightarrow{f} X_L \to X_M\to X_N \to Y_L \to Y_M\to Y_N \xrightarrow{f} X_L \cdots.
$$
Cutting the exact sequence at $f$ and using the additivity of length, we obtain 
$$
\ell(\image f)+\ell(X_M)+\ell(Y_L)+\ell(Y_N)=\ell(X_L)+\ell(X_N)+\ell(Y_M)+\ell(\image f).
$$
This gives $\theta(M)=\theta(L)+\theta(N)$.
Thus, it suffices to show that $\theta(R/\p)=0$ for every minimal prime $\p$ of $R$, since the claim then follows by a filtration argument.

Let $\p$ be a minimal prime of $R$.
Since $R/\p$ is a one-dimensional complete local domain, its normalization $S$ is a DVR and finitely generated over $R$.
As $S/(R/\p)$ has finite length, we have $\theta(S)=\theta(R/\p)$.
The $S$-module $\omega_R\otimes_R S$ has rank one, so we can write $\omega_R\otimes_R S\cong S/(a_1)\oplus \cdots \oplus S/(a_{t-1})\oplus S$ for some $a_1,\ldots,a_{t-1}\in S$ by the structure theorem over PID.
Hence $A_S$ is equivalent to a diagonal matrix $\operatorname{diag}(a_1,\ldots,a_{t-1},0)$.
For $b=\prod_{i=1}^{t-1}a_i$, $B_S$ is equivalent to a diagonal matrix $\operatorname{diag}(0,\ldots,0,b)$. One has
\begin{align*}
&X_S=\ker (\operatorname{diag}(a_1,\ldots,a_{t-1},0))/\image (\operatorname{diag}(0,\ldots,0,b))=S/bS, \\
&Y_S=\ker (\operatorname{diag}(0,\ldots,0,b))/\image(\operatorname{diag}(a_1,\ldots,a_{t-1},0))=\bigoplus_{i=1}^{t-1} S/a_iS.
\end{align*}
We have $\ell_S(X_S)=v(b)=\sum_{i=1}^{t-1} v(a_i)=\ell_S(Y_S)$, where $v$ denotes the valuation on $S$.
After canceling the factor corresponding to the degree of the field extension, we obtain $\ell_R(X_S)=\ell_R(Y_S)$, and thus $\theta(S)=0$.
In particular, we obtain $\theta(R)=0$, and hence $\ell(X_R)=\ell(Y_R)$. 
Since $Y_R=\ker B/\image A$ has finite length and it is a submodule of $\omega_R=\cok A$, we must have $Y_R=0$. 
It then follows from $\ell(X_R)=\ell(Y_R)$ that $X_R=0$. 
This means that \eqref{ABAB} is exact.
As $\Ext_R^1(\image A, \omega_R)=0$, the induced exact sequence $0\to \omega_R \to R^{\oplus t} \to \image A\to 0$ splits.
This contradicts the assumption that $R$ is not Gorenstein.
\end{proof}

The above result removes the assumption that the ring has type $2$ in \cite[Corollary 6.7]{LM20}, thereby generalizing the result to arbitrary generically Gorenstein rings.
Question \ref{growth ques canonical}(4) has an affirmative answer whenever there exists a lift to a generically Gorenstein ring as follows.

\begin{cor}\label{monomial0<1}
Let $k[x_1,\ldots,x_n]$ be a polynomial ring over $k$ and let $I$ be a monomial ideal such that $R=k[x_1,\ldots,x_n]/I$ is a non-Gorenstein Cohen--Macaulay ring. Then $\beta_1(\omega_R)>\beta_0(\omega_R)$.
\end{cor}

\begin{proof}
As in Proposition \ref{monomial}, by combining Remarks \ref{rmkques} and \ref{grloc}, we may assume that $I$ is squarefree. Since $R$ is reduced, we may reduce to the local case and apply Theorem \ref{genGor0<1}.
\end{proof}

As a consequence of Lemma \ref{samemonosameBetti} and Corollary \ref{monomial0<1}, we obtain the following result.

\begin{cor}
Let $(Q,\n,k)$ be a Gorenstein local ring, let $x_1,\ldots,x_n$ be a $Q$-regular sequence, and let $I$ be an ideal of $Q$ generated by monomials in $x_1,\ldots,x_n$.
If $R=Q/I$ is a non-Gorenstein Cohen--Macaulay ring, then $\beta_1(\omega_R)>\beta_0(\omega_R)$.
\end{cor}

An ideal is called a \textit{complete intersection} if it is generated by a regular sequence, and it is called \textit{generically a complete intersection} if it is a complete intersection locally at each of its associated prime ideals.
Since every Cohen--Macaulay local ring of codimension $2$ can be expressed, after completion, as a quotient of a regular local ring by a perfect ideal of height $2$, which is necessarily licci, the following corollary refines \cite[Proposition 2.7]{JL07}.

\begin{cor}\label{licci0<1}
Let $(Q,\n)$ be a Gorenstein local ring and let $I$ be in the linkage class of a generically complete intersection such that $R=Q/I$ is Cohen--Macaulay (e.g.\! $I$ is licci).
If $R$ is not Gorenstein, then $\beta_1(\omega_R)>\beta_0(\omega_R)$.
\end{cor}

\begin{proof}
Passing to the faithfully flat local extension $Q\to Q[x]_{\n Q[x]}$ with Gorenstein closed fiber, we may assume that the residue field of $Q$ is infinite.
It follows from \cite[Lemma 4.7]{JM14} that there exist a local ring $S$, a generically complete intersection ideal $J$ of $S$, and a regular sequence $x_1,\ldots,x_n$ on both $S$ and $T:=S/J$ such that $S/(x_1,\ldots,x_n)S\cong Q$ and $T/(x_1,\ldots,x_n)T\cong R$.
Since $Q$ is Gorenstein and $R$ is Cohen–Macaulay, $S$ and $T$ are Gorenstein and Cohen–Macaulay, respectively.
In particular, $T$ admits a canonical module $\omega_T$.
As $T$ is generically Gorenstein, we have $\beta_1(\omega_R)=\beta_1(\omega_T)>\beta_0(\omega_T)=\beta_0(\omega_R)$ by Remark \ref{rmkques} and Theorem \ref{genGor0<1}.
\end{proof}

Finally, we conclude the paper with a general result on the Betti numbers of finitely generated modules of finite injective dimension, which holds whenever Question \ref{growth ques canonical}(4) has an affirmative answer.

\begin{prop}\label{everyfid0<1}
Let $R$ be a Cohen--Macaulay local ring with canonical module $\omega_R$.
If $\beta_1(\omega_R)>\beta_0(\omega_R)$, then $\beta_1(M)>\beta_0(M)$ for every nonzero finitely generated $R$-module $M$ with finite injective dimension.
\end{prop}

\begin{proof}
Set $N=\Hom_R(\omega_R,M)$.
By the Foxby equivalence, $M\cong \omega_R\otimes_R N$ and $\Tor_i^R(\omega_R,N)=0$ for all $i>0$.
Hence the tensor product of minimal free resolutions of $\omega_R$ and $N$ yields a minimal free resolution of $M$.
By assumption, one has $\beta_0(M)=\beta_0(\omega_R)\beta_0(N)<\beta_1(\omega_R)\beta_0(N)+\beta_0(\omega_R)\beta_1(N)=\beta_1(M)$.
\end{proof}

\begin{ac}
The author was supported by JSPS Overseas Research Fellowships.
The author is grateful to Linquan Ma, the host researcher for the JSPS program, for suggesting various directions for research and for his valuable comments.
The author thanks ChatGPT for providing the content of Remark \ref{notopenle}. The same AI was also used to assist with English writing and mathematical discussions.
\end{ac}

\end{document}